\input phyzzx
\catcode`@=11 
\def\space@ver#1{\let\@sf=\empty \ifmmode #1\else \ifhmode
   \edef\@sf{\spacefactor=\the\spacefactor}\unskip${}#1$\relax\fi\fi}
\def\attach#1{\space@ver{\strut^{\mkern 2mu #1} }\@sf\ }
\newtoks\foottokens
\newbox\leftpage \newdimen\fullhsize \newdimen\hstitle
\newdimen\hsbody
\newif\ifreduce  \reducefalse
\def\almostshipout#1{\if L\lr \count2=1
      \global\setbox\leftpage=#1 \global\let\lr=R
  \else \count2=2
    \shipout\vbox{\special{dvitops: landscape}
      \hbox to\fullhsize{\box\leftpage\hfil#1}} \global\let\lr=L\fi}
\def\smallsize{\relax
\font\eightrm=cmr8 \font\eightbf=cmbx8 \font\eighti=cmmi8
\font\eightsy=cmsy8 \font\eightsl=cmsl8 \font\eightit=cmti8
\font\eightt=cmtt8
\def\eightpoint{\relax
\textfont0=\eightrm  \scriptfont0=\sixrm
\scriptscriptfont0=\sixrm
\def\rm{\fam0 \eightrm \f@ntkey=0}\relax
\textfont1=\eighti  \scriptfont1=\sixi
\scriptscriptfont1=\sixi
\def\oldstyle{\fam1 \eighti \f@ntkey=1}\relax
\textfont2=\eightsy  \scriptfont2=\sixsy
\scriptscriptfont2=\sixsy
\textfont3=\tenex  \scriptfont3=\tenex
\scriptscriptfont3=\tenex
\def\it{\fam\itfam \eightit \f@ntkey=4 }\textfont\itfam=\eightit
\def\sl{\fam\slfam \eightsl \f@ntkey=5 }\textfont\slfam=\eightsl
\def\bf{\fam\bffam \eightbf \f@ntkey=6 }\textfont\bffam=\eightbf
\scriptfont\bffam=\sixbf   \scriptscriptfont\bffam=\sixbf
\def\tt{\fam\ttfam \eightt \f@ntkey=7 }
\def\caps{\fam\cpfam \tencp \f@ntkey=8 }\textfont\cpfam=\tencp
\setbox\strutbox=\hbox{\vrule height 7.35pt depth 3.02pt width\z@}
\samef@nt}
\def\Eightpoint{\eightpoint \relax
  \ifsingl@\subspaces@t2:5;\else\subspaces@t3:5;\fi
  \ifdoubl@ \multiply\baselineskip by 5
            \divide\baselineskip by 4\fi }
\parindent=16.67pt
\itemsize=25pt
\thinmuskip=2.5mu
\medmuskip=3.33mu plus 1.67mu minus 3.33mu
\thickmuskip=4.17mu plus 4.17mu
\def\thinspace{\kern .13889em }
\def\negthinspace{\kern-.13889em }
\def\enspace{\kern.416667em }
\def\enskip{\hskip.416667em\relax}
\def\quad{\hskip.83333em\relax}
\def\qquad{\hskip1.66667em\relax}
\def\crr{\cropen{8.3333pt}}
\foottokens={\Eightpoint\singlespace}
\def\papersize{\SIZE\OFFSET\skip\footins=\bigskipamount}
\def\SIZE{\hsize=11.8truecm\vsize=17.5truecm}
\def\OFFSET{\voffset=-1.3truecm\hoffset=  .14truecm}
\message{STANDARD CERN-PREPRINT FORMAT}
\def\attach##1{\space@ver{\strut^{\mkern 1.6667mu ##1} }\@sf\ }
\def\PH@SR@V{\doubl@true\baselineskip=20.08pt plus .1667pt minus
.0833pt
             \parskip = 2.5pt plus 1.6667pt minus .8333pt }
\def\author##1{\vskip\frontpageskip\titlestyle{\tencp ##1}\nobreak}
\def\address##1{\par\kern 4.16667pt\titlestyle{\tenpoint\it ##1}}
\def\andaddress{\par\kern 4.16667pt \centerline{\sl and} \address}
\def\abstract{\vskip2\frontpageskip\centerline{\tenrm Abstract}
              \vskip\headskip }
\def\cases##1{\left\{\,\vcenter{\Tenpoint\m@th
    \ialign{$####\hfil$&\quad####\hfil\crcr##1\crcr}}\right.}
\def\matrix##1{\,\vcenter{\Tenpoint\m@th
    \ialign{\hfil$####$\hfil&&\quad\hfil$####$\hfil\crcr
      \mathstrut\crcr\noalign{\kern-\baselineskip}
     ##1\crcr\mathstrut\crcr\noalign{\kern-\baselineskip}}}\,}
\Tenpoint
}
\def\Smallsize{\smallsize\reducetrue
\let\lr=L
\hstitle=8truein\hsbody=4.75truein\fullhsize=24.6truecm\hsize=\hsbody
\output={
  \almostshipout{\leftline{\vbox{\makeheadline
  \pagebody\makefootline}}}\advancepageno
     }
\special{dvitops: landscape}
\def\makeheadline{
\iffrontpage\line{\the\headline}
             \else\vskip .0truecm\line{\the\headline}\vskip .5truecm
\fi}
\def\makefootline{\iffrontpage\vskip  0.truecm\line{\the\footline}
               \vskip -.15truecm\line{\the\date\hfil}
              \else\line{\the\footline}\fi}
\paperheadline={
\iffrontpage\hfil
               \else
               \tenrm\hss $-$\ \folio\ $-$\hss\fi    }
\paperstyle}
%
%
%
%
%
%
%
%
%
\newcount\referencecount     \referencecount=0
\newif\ifreferenceopen       \newwrite\referencewrite
\newtoks\rw@toks
\def\NPrefmark#1{\attach{\scriptscriptstyle [ #1 ] }}
\let\PRrefmark=\attach
\def\refmark#1{\relax\ifPhysRev\PRrefmark{#1}\else\NPrefmark{#1}\fi}
\def\refend{\refmark{\number\referencecount}}
\newcount\lastrefsbegincount \lastrefsbegincount=0
\def\refsend{\refmark{\count255=\referencecount
   \advance\count255 by-\lastrefsbegincount
   \ifcase\count255 \number\referencecount
   \or \number\lastrefsbegincount,\number\referencecount
   \else \number\lastrefsbegincount-\number\referencecount \fi}}
\def\refch@ck{\chardef\rw@write=\referencewrite
   \ifreferenceopen \else \referenceopentrue
   \immediate\openout\referencewrite=referenc.texauxil \fi}
%
{\catcode`\^^M=\active 
  \gdef\obeyendofline{\catcode`\^^M\active \let^^M\ }}%
%
{\catcode`\^^M=\active 
  \gdef\ignoreendofline{\catcode`\^^M=5}}
{\obeyendofline\gdef\rw@start#1{\def\t@st{#1} \ifx\t@st\blankend%
\endgroup \@sf \relax \else \ifx\t@st\bl@nkend \endgroup \@sf \relax%
\else \rw@begin#1
\backtotext
\fi \fi } }
{\obeyendofline\gdef\rw@begin#1
{\def\n@xt{#1}\rw@toks={#1}\relax%
\rw@next}}
\def\blankend{}
{\obeylines\gdef\bl@nkend{
}}
\newif\iffirstrefline  \firstreflinetrue
\def\rwr@teswitch{\ifx\n@xt\blankend \let\n@xt=\rw@begin %
 \else\iffirstrefline \global\firstreflinefalse%
\immediate\write\rw@write{\noexpand\obeyendofline \the\rw@toks}%
\let\n@xt=\rw@begin%
      \else\ifx\n@xt\rw@@d \def\n@xt{\immediate\write\rw@write{%
        \noexpand\ignoreendofline}\endgroup \@sf}%
             \else \immediate\write\rw@write{\the\rw@toks}%
             \let\n@xt=\rw@begin\fi\fi \fi}
\def\rw@next{\rwr@teswitch\n@xt}
\def\rw@@d{\backtotext} \let\rw@end=\relax
\let\backtotext=\relax

\newdimen\refindent     \refindent=30pt
\def\refitem#1{\par \hangafter=0 \hangindent=\refindent
\Textindent{#1}}
\def\REFNUM#1{\space@ver{}\refch@ck \firstreflinetrue%
 \global\advance\referencecount by 1 \xdef#1{\the\referencecount}}
\def\refnum#1{\space@ver{}\refch@ck \firstreflinetrue%
 \global\advance\referencecount by 1
\xdef#1{\the\referencecount}\refend}

\def\REF#1{\REFNUM#1%
 \immediate\write\referencewrite{%
 \noexpand\refitem{#1.}}%
\begingroup\obeyendofline\rw@start}
\def\ref{\refnum\?%
 \immediate\write\referencewrite{\noexpand\refitem{\?.}}%
\begingroup\obeyendofline\rw@start}
\def\Ref#1{\refnum#1%
 \immediate\write\referencewrite{\noexpand\refitem{#1.}}%
\begingroup\obeyendofline\rw@start}
\def\REFS#1{\REFNUM#1\global\lastrefsbegincount=\referencecount
\immediate\write\referencewrite{\noexpand\refitem{#1.}}%
\begingroup\obeyendofline\rw@start}
\newif\ifmref  
\newif\iffref  
\def\xrefsend{\xrefmark{\count255=\referencecount
\advance\count255 by-\lastrefsbegincount
\ifcase\count255 \number\referencecount
\or \number\lastrefsbegincount,\number\referencecount
\else \number\lastrefsbegincount-\number\referencecount \fi}}
\def\xrefsdub{\xrefmark{\count255=\referencecount
\advance\count255 by-\lastrefsbegincount
\ifcase\count255 \number\referencecount
\or \number\lastrefsbegincount,\number\referencecount
\else \number\lastrefsbegincount,\number\referencecount \fi}}
\def\xREFNUM#1{\space@ver{}\refch@ck\firstreflinetrue%
\global\advance\referencecount by 1
\xdef#1{\xrefend}}
\def\xrefend{\xrefmark{\number\referencecount}}
\def\xrefmark#1{[{#1}]}
\def\xRef#1{\xREFNUM#1\immediate\write\referencewrite%
{\noexpand\refitem{#1 }}\begingroup\obeyendofline\rw@start}%
\def\xREFS#1{\xREFNUM#1\global\lastrefsbegincount=\referencecount%
\immediate\write\referencewrite{\noexpand\refitem{#1 }}%
\begingroup\obeyendofline\rw@start}
\def\rrr#1#2{\relax\ifmref{\iffref\xREFS#1{#2}%
\else\xRef#1{#2}\fi}\else\xRef#1{#2}\xrefend\fi}
\referencecount=0
\def\refbreak{\hfil\penalty200\hfilneg}
\def\paperstyle{\papers}
\paperstyle   
%
%
%
\def\slacpub{\afterassignment\slacp@b\toks@}
\def\slacp@b{\edef\n@xt{\Pubnum={\the\toks@}}\n@xt}
\let\pubnum=\slacpub
\expandafter\ifx\csname eightrm\endcsname\relax
    \let\eightrm=\ninerm \let\eightbf=\ninebf \fi

\font\seventeencp=cmcsc10 scaled\magstep3

\newif\ifCONF \CONFfalse
\newif\ifBREAK \BREAKfalse
\newif\ifsectionskip \sectionskiptrue

%
%
%
%
\def\NuclPhysProc{
\let\lr=L
\hstitle=8truein\hsbody=4.75truein\fullhsize=21.5truecm\hsize=\hsbody
\hstitle=8truein\hsbody=4.75truein\fullhsize=20.7truecm\hsize=\hsbody
\output={
  \almostshipout{\leftline{\vbox{\makeheadline
  \pagebody\makefootline}}}\advancepageno
     }
\def\papersize{\SIZE\OFFSET\skip\footins=\bigskipamount}
\def\SIZE{\hsize=10.0truecm\vsize=27.0truecm}
\def\OFFSET{\voffset=-1.4truecm\hoffset=-2.40truecm}
\message{NUCLEAR PHYSICS PROCEEDINGS FORMAT}
\def\makeheadline{
\iffrontpage\line{\the\headline}
             \else\vskip .0truecm\line{\the\headline}\vskip .5truecm
\fi}
\def\makefootline{\iffrontpage\vskip  0.truecm\line{\the\footline}
               \vskip -.15truecm\line{\the\date\hfil}
              \else\line{\the\footline}\fi}
\paperheadline={\hfil}
\paperstyle}
%
%
%
%

%
%
%
%

%
%
%
%
\def\ReprintVolume{\smallsize
\def\papersize{\hsize=18.0truecm\vsize=23.1truecm\voffset -.73truecm
    \hoffset -.65truecm\skip\footins=\bigskipamount
    \normaldisplayskip= 20pt plus 5pt minus 10pt}
\message{REPRINT VOLUME FORMAT}
\paperstyle\baselineskip=.425truecm\parskip=0truecm
\def\makeheadline{
\iffrontpage\line{\the\headline}
             \else\vskip .0truecm\line{\the\headline}\vskip .5truecm
\fi}
\def\makefootline{\iffrontpage\vskip  0.truecm\line{\the\footline}
               \vskip -.15truecm\line{\the\date\hfil}
              \else\line{\the\footline}\fi}
\paperheadline={
\iffrontpage\hfil
               \else
               \tenrm\hss $-$\ \folio\ $-$\hss\fi    }
\def\sectionfont{\bf}    }
%
%
%
%
\def\SIZE{\hsize=15.73truecm\vsize=23.11truecm}
\def\OFFSET{\voffset=0.0truecm\hoffset=0.truecm}
\message{DEFAULT FORMAT}
\def\papersize{\SIZE\OFFSET\skip\footins=\bigskipamount
\normaldisplayskip= 35pt plus 3pt minus 7pt}
\Pubnum={\rm \the\pubnum }
\def\title#1{\vskip\frontpageskip\vskip .50truein
     \titlestyle{\seventeencp #1} \vskip\headskip\vskip\frontpageskip
     \vskip .2truein}
\def\author#1{\vskip .27truein\titlestyle{#1}\nobreak}

\def\p@bblock{\begingroup \tabskip=\hsize minus \hsize
   \baselineskip=1.5\ht\strutbox \topspace+2\baselineskip
   \halign to\hsize{\strut ##\hfil\tabskip=0pt\crcr
  \the \Pubnum\cr}\endgroup}
\def\makefootline{\iffrontpage\vskip .27truein\line{\the\footline}
                 \vskip -.1truein
              \else\line{\the\footline}\fi}
\paperfootline={\iffrontpage\message{FOOTLINE}
\hfil\else\hfil\fi}

\def\abstract{\vskip2\frontpageskip\centerline{\twelvebf Abstract}
              \vskip\headskip }

\paperheadline={
\iffrontpage\hfil
               \else
               \twelverm\hss $-$\ \folio\ $-$\hss\fi}
%
%
\def\nup#1({\refbreak\ Nucl.\ Phys.\ $\underline {B#1}$\ (}
\def\plt#1({\refbreak\ Phys.\ Lett.\ $\underline  {#1}$\ (}
\def\cmp#1({\refbreak\ Commun.\ Math.\ Phys.\ $\underline  {#1}$\ (}
\def\prp#1({\refbreak\ Physics\ Reports\ $\underline  {#1}$\ (}
\def\prl#1({\refbreak\ Phys.\ Rev.\ Lett.\ $\underline  {#1}$\ (}
\def\prv#1({\refbreak\ Phys.\ Rev. $\underline  {D#1}$\ (}
\def\und#1({            $\underline  {#1}$\ (}
%
%

\def\rB{\hfil\penalty1000\hfilneg}
%
%
\hyphenation{sym-met-ric anti-sym-me-tric re-pa-ra-me-tri-za-tion
Lo-rentz-ian a-no-ma-ly di-men-sio-nal two-di-men-sio-nal}
%
%
%
%

\def\coeff#1#2{{\textstyle { #1 \over #2}}\displaystyle}
\def\boxit#1{\vbox{\hrule\hbox{\vrule\kern3pt
\vbox{\kern3pt#1\kern3pt}\kern3pt\vrule}\hrule}}
\message{ by V.K, W.L and A.S}
\catcode`@=12
\paperstyle
\input tables

\def\chi {\X}

\paperstyle

\def\mod{{\rm ~mod~}}

\def\half{\coeff12}

\def\Zbf{{\bf Z}}
\def\CFT{{\rm C}}

\def\X{{\cal X}}

\catcode`@=11
\def\ninef@nts{\relax
    \textfont0=\ninerm          \scriptfont0=\sixrm
      \scriptscriptfont0=\sixrm
    \textfont1=\ninei           \scriptfont1=\sixi
      \scriptscriptfont1=\sixi
    \textfont2=\ninesy          \scriptfont2=\sixsy
      \scriptscriptfont2=\sixsy
    \textfont3=\tenex          \scriptfont3=\tenex
      \scriptscriptfont3=\tenex
    \textfont\itfam=\nineit     \scriptfont\itfam=\seveni  
\sevenit
    \textfont\slfam=\ninesl     \scriptfont\slfam=\sixrm 
\sevensl
    \textfont\bffam=\ninebf     \scriptfont\bffam=\sixbf
      \scriptscriptfont\bffam=\sixbf
    \textfont\ttfam=\tentt
    \textfont\cpfam=\tencp }
\def\ninepoint{\ninef@nts \samef@nt \b@gheight=9pt \setstr@t }
\newif\ifnin@  \nin@false
\def\Tenpoint{\tenpoint\twelv@false\nin@false\spaces@t}
\def\Twelvepoint{\twelvepoint\twelv@true\nin@false\spaces@t}
\def\Ninepoint{\ninepoint\twelv@false\nin@true\spaces@t}
\def\spaces@t{\rel@x
      \iftwelv@ \ifsingl@\subspaces@t3:4;\else\subspaces@t1:1;\fi
       \else \ifsingl@\subspaces@t3:5;\else\subspaces@t4:5;\fi \fi
      \ifdoubl@ \multiply\baselineskip by 5
         \divide\baselineskip by 4 \fi
       \ifnin@ \ifsingl@\subspaces@t3:8;\else\subspaces@t4:7;\fi \fi
}
\def\Vfootnote#1{\insert\footins\bgroup
   \interlinepenalty=\interfootnotelinepenalty \floatingpenalty=20000
   \singl@true\doubl@false \iftwelv@ \Tenpoint
   \else \Ninepoint \fi
   \splittopskip=\ht\strutbox \boxmaxdepth=\dp\strutbox
   \leftskip=\footindent \rightskip=\z@skip
   \parindent=0.5\footindent \parfillskip=0pt plus 1fil
   \spaceskip=\z@skip \xspaceskip=\z@skip \footnotespecial
   \Textindent{#1}\footstrut\futurelet\next\fo@t}

\def\small#1{\vskip .3truecm\footnoterule\nobreak
\Ninepoint\parindent=2pc\sl
\hang #1 \vskip .3truecm\nobreak\footnoterule}

\def\small#1{}

  

\def\GepF{\rrr\GepF{
D.~Gepner, \plt B222 (1989) 207.}}
\def\ScYe{\rrr\ScYe{A.N.~Schellekens and S.~Yankielowicz,
\nup334 (1990) 67.}}
\def\ScYg{\rrr\ScYg{A.N.~Schellekens and S.~Yankielowicz,
Int.J.Mod.Phys.\und{A5} (1990) 2903.}}
\def\MoSc{\rrr\MoSc{
G.~Moore   and N.~Seiberg,
\plt B220 (1989) 422.}}
\def\LVW {\rrr\LVW {
W.~Lerche, C.~Vafa and N.~Warner,
\nup 324 (1989) 427.}}
\def\BeBU{\rrr\BeBU {B.~Gato-Rivera and A.N.~Schellekens,
Commun.~Math.~Phys.\und{145} (1992) 85;\rB
M.~Kreuzer and A.N.~Schellekens,
\nup411 (1994) 97.}}
\def\FSSc{\rrr\FSSc{ J. Fuchs, A.N. Schellekens and C. Schweigert,
\cmp 180 (1996) 39.}}
\def\DuJo{\rrr\DuJo{D. Dunbar and K. Joshi, Mod.Phys.Lett.A8:2803-2814,1993.}}
\def\FSSs{\rrr\FSSs{
J. Fuchs, A.N. Schellekens and C. Schweigert,
\nup 473 (1996) 323.}}
\def\FSSt{\rrr\FSSt{
J. Fuchs, A.N. Schellekens and C. Schweigert,
\nup 461 (1996) 371.}}
\def\MoSeCMP{\rrr\MoSeCMP{
G. Moore and N. Seiberg,  Commun.Math.Phys. 123 (1989) 177.}}
\def\Bant{\rrr\Bant{
P. Bantay,  Int. J. Mod. Phys. A13 (1998) 175.}}
\def\FuSc{\rrr\FuSc{
J.~Fuchs and C.~Schweigert,  Phys. Lett. B414 (1997) 251;\rB
Phys. Lett. B447 (1999) 266;\rB {\it ~~Symmetry Breaking Boundaries. 1.  
General Theory.}\rB~~Preprint CERN-TH-99-35 (hep-th/9902132)}}
\def\PSS{\rrr\PSS{
G.~Pradisi, A.~Sagnotti and Ya.S.~Stanev, Phys Lett B354 (1995) 279;
Phys. Lett B356 (1995) 230;
 Phys. Lett. B381 (1996) 97.}}

\pubnum={{}}
\rightline{NIKHEF/99-016}
\rightline{hep-th/9905153}
\rightline{May 1999}
\date{May 1999}
\pubtype{CRAP}
\titlepage
\message{TITLE}

\title{\fourteenbf Fixed point resolution in extended WZW-models }
\author{A.N. Schellekens\foot{t58@nikhef.nl}}
\line{\hfil
 {\it NIKHEF}
 \hfil}
\line{\hfil \it NIKHEF, P.O. Box 41882, 1009 DB Amsterdam,
The Netherlands  \hfil}
\bigskip

\abstract \noindent
A formula is derived for the fixed point resolution matrices
of simple current extended WZW-models and coset conformal
field theories. Unlike the analogous matrices for
unextended WZW-models, these matrices are in general not
symmetric, and they may have field-dependent twists. They
thus provide non-trivial realizations of the general
conditions presented in earlier work with Fuchs and Schweigert.

\baselineskip= 15.0pt plus .2pt minus .1pt

\chapter{Introduction}

Despite a large amount of work on conformal field theory during
the last fifteen years, there are still no efficient procedures
for computing many quantities of interest. Examples
of such quantities are the spectrum, the (Virasoro or extended)
characters, the modular transformation properties of these characters,
the fusion rules, fusing and braiding matrices, and correlation
functions on arbitrary Riemann surfaces.
While the difficulties in computing the latter quantities are
well-known, even the determination of the spectrum can involve
unexpected difficulties that are often overlooked in the literature.
Indeed, for
a large class of rational CFT's that are very explicitly defined,
namely the class of unitary coset conformal field theories, it is not even
known {\it in general} how to compute the spectrum (although in certain
simple cases no problems occur). Even determining the correct
{\it number} of primaries in such theories -- not to mention the
conformal weights -- is non-trivial in some classes of
coset theories and extended WZW-models, and
was only solved a few years ago in \FSSs.
The determination of the characters -- and hence the spectrum --
is only fully known for the special case of diagonal coset theories
\FSSt.

The problems alluded to in the end of the previous paragraph
all have a common origin, that of resolving fixed points in
simple current extensions. This problem
arises in many coset
CFT's (as was first noted in \LVW\ and \MoSc)
due to the procedure known as ``field identification"
\GepF, which
is formally analogous to a simple current extension of
the chiral algebra \ScYe.
In general -- whether applied to the construction of
extended chiral algebras or to field identification
-- simple currents may
have fixed points which must be resolved in order to arrive
at a correct description of the theory. This fixed point
resolution procedure requires data in addition to
the characters (or branching functions) and modular matrix $S$
of the unextended theory. This information is then used
to modify the coset branching functions to obtain the characters,
and to modify the ``naive" matrix $S$ to obtain the correct one.

It is the latter aspect of fixed point resolution that interests
us here. The additional information needed to resolve
the fixed points are the ``fixed point resolution matrices"
$S^J$, which exist for any simple current $J$ that has
fixed points. In \FSSs\ these matrices were obtained
for unextended WZW-models (``A-invariants").
They turned out to be simply
related to ordinary modular transformation matrices $S$
of an ``orbit Lie algebra" \FSSc\ related to the WZW-model.
In general
the matrices $S^J$ have been conjectured \FSSs\ to be related to the modular
transformation matrices of holomorphic one-point blocks \MoSeCMP.
Strong evidence for this conjecture was provided in \Bant.

In this paper we will give a formula for
the fixed point resolution matrices of
simple current extended rational conformal field theories. These
matrices are expressed in terms of the
corresponding matrices of the unextended CFT. Although
{\it a priori} nothing
restricts us to WZW-models, this is the only class for which
the fixed point resolution matrices are already known, and hence
the obvious starting point.
The main application of our result is therefore to
simple current extended
WZW-models (also known as
D-invariants) and field identification in coset theories.
This then generalizes the results of \FSSs\ to
all coset theories, except a
small class described in \DuJo, which has ``exceptional" field
identifications.
Since
most rational CFT's can be described as coset CFT's, this is a step
in the direction of generalizing the results of \FSSs\ to all
rational CFT's, and allows us to examine whether, and how,
the general case differs from the special case of WZW-models.

In \FSSs\ a set of conditions for fixed point
resolution matrices in generic rational CFT's was presented.
The fixed point resolution matrices of unextended WZW-models
were found to satisfy
additional constraints,
due to the fact that they are modular transformations of an
orbit Lie algebra.
We find that in
extended WZW-models
most of these additional constraints do not hold,
and only the necessary conditions are satisfied.
In particular the fixed point resolution matrices do not
appear to be related to any kind of orbit Lie algebra. This
concept apparently does not generalize to conformal field theory.

The fixed point resolution matrices $S^J$
are needed to compute the modular transformation matrix $S$
of simple current extended CFT's. Using the Verlinde formula
one may then compute the fusion rules. The results presented
here allow the computation of $S$ if a simple current extended
CFT is extended once more by another set of simple currents.
The result of this procedure is an enlargement
of the original chiral algebra in two steps. Obviously,
the final result
can also be obtained by making the entire extension in just
one step, and in that case the formula of \FSSs\ already gives
the matrix $S$ of the twice-extended theory in terms of
the fixed point resolution matrices of the unextended theory.
Hence the fixed point resolution matrices of the once-extended
theory are not strictly needed for the purpose of computing $S$
(although they are in some cases convenient for practical reasons).
However, it has recently become clear that the fixed point
resolution matrices appear also in other contexts \FuSc, namely
in the description of boundary conditions for open strings
built as ``descendants" of certain conformal field theories \PSS.
For these applications we need to know these matrices explicitly.

Since the general formula of \FSSs\ covers the
first extension, the second extension as well as the
full extension, it is clear that the information we
are looking for can be derived from the results of that paper.
Indeed, this fact played an important r\^ole in obtaining
the general formula and the conditions
governing $S^J$, and some cases of successive extensions
are discussed in \FSSs.
It turns out, however, that there
are several subtleties in the general case. After treating
these carefully, we finally arrive at a formula for $S^J$.
The fixed point resolution matrices $S^J$ may be viewed
as a generalization of the modular  transformation
matrix $S\equiv S^1$. It is then not surprising that
the formula we obtain for $S^J$
is quite similar to the one obtained for $S$ in \FSSs.
The main difference consists of some additional
phase factors that vanish in the special case of $S^1$, but
are essential for consistency of $S^J$.

This paper is organized as follows. In the next section
we summarize the results of \FSSs, with a minor
improvement in the formulation of the conditions for $S^J$.
The conditions involve two new quantities, the ``simple
current twist" $F(a,K,J)$ and the ``conjugation twist" $G(a,K,J)$,
where $a$ is a primary field fixed by the simple currents $K$ and $J$.
These quantities are equal to each other for unextended WZW-models
(and independent of $a$). We give an argument for the validity of this
equality in all CFT's with real twists (which so far are the
only ones known). The
details are given in Appendix A.
On the other hand, $a$-independence turns out
{\it not} to be a general feature, as is shown by means of an example in
section 5.
In section 3 we give the precise formulation of the
problem we want to solve, and we also discuss the
choice of orbit representatives needed in the derivation.
The existence proof for our choice of representatives is presented
in Appendix B.
The derivation of the main result is presented in section 4. It involves
two main ingredients: the proper choice of representatives, and
the splitting of discrete group characters. The mathematical
background for this splitting is provided in Appendix C. Both of
these ingredients lead to crucial phase factors in the
main result, formula 4.4.
In section 5 we investigate some of the properties of the
matrices $S^J$.

\chapter{Fixed point resolution}

In this section we summarize the results of \FSSs. The
reader is assumed to be familiar with simple currents and
fixed points (see \ScYg\ for a review).

Consider a rational CFT ${\CFT}$ with a set of
mutually local, integer spin simple currents forming
a subgroup ${\cal H}$ of the full center ${\cal G}$.
We extend the chiral algebra by the currents in ${\cal H}$,
and we denote the new CFT as ${\CFT}^{\cal H}$.
The fields in the new theory correspond to
${\cal H}$-orbits of fields\foot{The term ``field" is used
throughout this paper as an abbreviation for ``primary field"}of
${\CFT}$ that are local with
respect to ${\cal H}$, but in general this is not a one-to-one
correspondence. For each field $a$ of ${\CFT}$ we
define the stabilizer ${\cal S}_a$ as the set of currents $J$
in ${\cal H}$ that fixes $a$: $Ja=a$. Furthermore we define
a certain subgroup of ${\cal S}_a$, the untwisted stabilizer
${\cal U}_a$. The details of this definition follow below.
Then the fields of ${\CFT}^{\cal H}$ are labelled by a pair
$(a,i)$, where $a$ is a representative of an
${\cal H}$-local orbit, and $i$ is a character label
of ${\cal U}_a$. The modular transformation matrix
of $\CFT^{\cal H}$ is given by \FSSs

$$ S_{(a,i)(b,j)}= { | {\cal H} | \over
\sqrt{ | {\cal S}_a | | {\cal U}_a |
| {\cal S}_b | | {\cal U}_b |}} \sum_{J}
\Psi_{i}^J S^J_{ab} (\Psi_{j}^J)^* \ , \eqn\FSSform $$
where the sum is over the intersection of the two
untwisted stabilizers. The characters should carry an extra label
 $a$ or $b$
to indicate to which untwisted stabilizer they belong, but this
label is always already implied by the fixed point resolution labels,
and will therefore be suppressed as much as possible.

The matrices $S^J_{ab}$ appearing in this formula must
satisfy the following properties
\item\bullet $\{1\}\ \ \ S^J_{ab}=0$ if $Ja\not=a$ or $Jb\not=b$
\item\bullet $\{2\}\ \ \ S^J(S^J)^{\dagger}=1$\item\bullet $\{3\}\ \ \  
(S^JT^J)^3=(S^J)^{2}$
\item\bullet $\{4\}\ \ \ S^J_{Ka,b}=F(a,K,J)e^{2\pi i Q_K(b)}S^J_{ab}$
\item\bullet $\{4a\}\ \ \ F(a,K,J_1)F(a,K,J_2)=F(a,K,J_1J_2)$
\item\bullet $\{5\}\ \ \ (S^J)^{2}=\eta^J C^J$
\item\bullet $\{5a\}\ \ \
 \eta^J_{ab}=\eta^J_a \delta_{ab}$
\item\bullet $\{5b\}\ \ \  \eta_a^J\eta_a^K=G(a,K,J)\eta_a^{JK}\ , \ \ \ \
G(a,K,J) = 1\ \  \hbox{if}\ F(a,K,J) = 1$
\item\bullet $\{5c\}\ \ \ \eta^JC^J=C^J(\eta^J)^*$
\item\bullet $\{6\} \ \ \ S^J_{ab}=S^{J^{-1}}_{ba} $

Here properties $\{2\}-\{6\}$
are defined on the subspace of fields where $S^J$ is non-zero,
and $T^J$, $C^J$ are
the restrictions of $T$ and $C$ to that subspace; as usual $T$
is the generator of the $\tau \to \tau+1$ transformation of
the modular group, $S$ is the generator of $\tau \to - {1\over\tau}$,
and $C$ is the charge conjugation matrix. These  three matrices
satisfy $(ST)^3=S^2=C$. The quantity $Q_J(a)$ is the monodromy
charge of the field $a$ w.r.t. the current $J$.
The untwisted stabilizer ${\cal U}_a$ is defined as
$$ {\cal U}_a:=\{ J \in {\cal S}_a\ |\ F(a,K,J)=1\ \hbox{for all}\
K \in {\cal S}_a \} $$.

In comparison to \FSSs\ the condition $\{5b\}$ has been formulated
differently. In \FSSs\ it
was given as
$\eta_a^J\eta_a^K=\eta_a^{JK}$ if $J,K \in {\cal U}_a$.
It is easy to see that both forms of the condition are
equivalent, but the form chosen here has the advantage that no
explicit mention is made of (untwisted) stabilizers in the conditions
for $S^J$.
Indeed,
the matrices $S^J$ are a property of the CFT ${\CFT}$ and do
not depend on the extension of the chiral algebra that one considers,
and hence in particular not on ${\cal S}_a$ and ${\cal U}_a$.

One can derive from these conditions that the simple
current twists $F(a,K,J)$
are phases satisfying
$$ F(a,K_1,J)F(a,K_2,J)=F(a,K_1K_2,J) \eqn\faconeprod $$
and
$$ F(a,K,J)=F(a,J,K)^*\ , \eqn\fsym $$
provided that $Ka=a$
(note that the definition $\{4\}$ of $F$
requires $J$, but not $K$, to fix $a$)\rlap.
\foot{If $Ka\not=a$ the twists depend on a
 fixed point
labelling convention: one may label the resolved fixed points
of $a$ differently than those of $Ka$.}
Furthermore one can show that $F(a,J,J)=1$ if $J$
has integer spin, and $F(a,J,J)=-1$ if $J$ has half-integer spin.

In the special case  of unextended
WZW-models the matrices $S^J$ were obtained in \FSSs.
They have some additional
properties not required for fixed point resolution, and in particular
they are themselves modular transformation matrices of a set of
characters of (in some cases twisted) affine Lie algebras. Therefore
they are symmetric, and $S^J=S^{J'}$ if $J$ and $J'$ generate the same
cyclic subgroup   (hence
$S^J=S^{J^{-1}}$). Furthermore
 $\eta^J_a$ and $F(a,K,J)$ are independent of $a$.
For unextended WZW-models the twists satisfy the empirical relation
$$ G(a,K,J)=F(a,K,J)\ . \eqn\GisF $$
which implies $\{5b\}$, but is stronger.

Note that \GisF\ implies that $F$ is real, \ie\ $F=\pm1$, since
$G$ is symmetric in $J$ and $K$, and $F$ is symmetric only if it
is real. Under a very mild additional assumption the converse
is also true, \ie\ reality of $F$ implies \GisF. The additional
assumption is
\item\dash The matrices $S^J$ of tensor product CFT's are the
tensor products of the matrices $S^J$ of the factors.

\noindent
The proof of the converse
statement (\ie\ reality of $F$ implies $F=G$)
is given in appendix A. In all known
cases $F$ is real, and \GisF\ holds. It would be
interesting to find a fundamental reason why either property
should hold in general, but at least this argument shows that they
are related.

\chapter{Formulation of the problem}

\section{Definitions}

Consider a CFT $\CFT$ with center ${\cal G}$.
The fields of the CFT will be denoted by $a$.
For each field $a$ there is a subgroup of ${\cal G}$ that
fixes $a$. This group will be called the {\it full stabilizer}
of $a$, and will be denoted as ${\cal T}_a$. It depends in fact
only on the ${\cal G}$-orbit of $a$, not on orbit members
individually. It consists of currents of integer or
half-integer spin.

For each current $J$ in ${\cal G}$ we assume the existence of a
matrix $S^J_{ab}$ satisfying the conditions of \FSSs,
stated in the previous section.

We will consider here new conformal field theories obtained
by extending the chiral algebra of $\CFT$ by means of a subgroup
${\cal H}$ of ${\cal G}$ consisting of mutually local
integer spin currents.
 The new theory will be called $\CFT^{\cal H}$. In general,
all quantities in the extended theory will be denoted by
superscripts ${\cal H}$. For example, $a^{\cal H}$ denotes
an ${\cal H}$ primary.
The center of the new theory will be denoted as
${\cal G}^{\cal H}$. It contains the subgroup
of ${\cal G}$ of currents local w.r.t. ${\cal H}$, organized
into ${\cal H}$-orbits\rlap.\foot{In rare cases ${\cal G}^{\cal H}$
may be larger due to extra currents originating from resolved
fixed points, but we will ignore such currents here.}
To determine the spectrum and the modular transformation
matrix of $\CFT^{\cal H}$ one needs the stabilizer of each field. This
is the set of currents in ${\cal H}$ that fixes $a$. Hence
${\cal S}_a = {\cal T}_a \cap {\cal H}$.
This stabilizer should not be confused with the full
stabilizer ${\cal T}^{\cal H}_{a^{\cal H}}$ of the field $a^{\cal H}$
of $\CFT^{\cal H}$, consisting of the currents in  ${\cal G}^{\cal H}$
that fix $a^{\cal H}$.
Another important group is the untwisted stabilizer, whose
characters labels the resolved fixed points.
More precisely,
the fields $a^{\cal H}$ are
labelled by a set $(a,i)$, where $i$ labels a
representation of ${\cal U}_a $
and $a$ is a representative of a ${\cal H}$-local ${\cal H}$-orbit.

Our goal is the computation of the
matrices $S^{J^{\cal H}}_{a^{\cal H}b^{\cal H}}$ and
the determination of the twists
$F^{\cal H}(a^{\cal H},K^{\cal H},J^{\cal H})$.

\section{Choice of representatives}

Consider the action of the simple current
$J^{\cal H}$ on a resolved fixed point field
$a^{\cal H}=(a,i)$.
The  elements of the corresponding coset are
$\{ J, J h_1, \ldots J h_n \}$, where $J$ is some ${\cal H}$-local element
of ${\cal G}$ and
$h_i$ denote the elements of ${\cal H}$.
The current $J$ is the representative chosen to
identify the orbit. Suppose that in the
${\cal H}$-{\it extended} theory
the current $J^{\cal H}$ fixes $a^{\cal H}$.
Then $J$ must
map $a$ to another member of the ${\cal H}$ orbit of $a$, \ie\
$J a = h_i a$ for some $h_i \in {\cal H}$.
 Then $J h_i^{-1}$ fixes $a$.
Hence we can always find a representative that fixes $a$, and
we will denote this representative as $X_a(J)$, $X_a(J) \in J{\cal H}$.
In general $X_a$ depends on $a$, and there is no universal choice of
representatives that avoids such a dependence for all currents.
The complete set
of $\CFT$-currents
in the coset element $J^{\cal H}$ that fixes $a$ is of the form
$$ X_a(J) {\cal S}_a \ .$$

The existence of such a representative $X_a$ is necessary, but
not sufficient to conclude that $J^{\cal H}$ fixes $a^{\cal H}$.
It only ensures that $J^{\cal H}$ fixes the first of the two
labels $(a,i)$. The action of $J^{\cal H}$ on the second label
is trivial if and only if
$$  F(a,X_a(J),K)=1 \hbox{~for all~} K \in {\cal U}_a \eqn\FPcondtwo $$
This condition is the same for all choices $X_a(J) {\cal S}_a$,
because (using \faconeprod\ and the definition of the
untwisted stabilizer)
$$ F(a,X_a h,K) = F(a,X_a ,K)F(a,h,K)=F(a,X_a ,K)\ , $$
for all $h \in {\cal S}_a$.

If \FPcondtwo\ holds it is always possible to make a choice
$R_a(J)$
out of the set of representatives $X^a(J) {\cal S}_a$
such that
$$  F(a,R_a(J),K)=1 \hbox{~for all~} K \in {\cal S}_a \eqn\FPcondthree $$
This choice $R_a$ is unique up to multiplication by elements
of ${\cal U}_a$ and is essential in the derivation
in the next section. The proof that
such a choice is always possible is given in appendix B.

\chapter{Derivation of the main result}

The FCFT-matrix for a current $J^{\cal H}$ can
in principle be computed from the
general formula. Instead of extending $\CFT^{\cal H}$ we may
extend $\CFT$ by a set of simple currents ${\cal M} \subset {\cal G}$
in such a way that ${\cal M}$ contains all
${\cal G}$-elements of the form  $J^n {\cal H}$ (where $J$ is a
representative of the coset corresponding to  $J^{\cal H}$).
Eqn. \FSSform\ is in this case
$$ S_{(a,\alpha)(b,\beta)}= { | {\cal M} | \over
\sqrt{ | \hat{\cal S}_a | | \hat{\cal U}_a |
| \hat{\cal S}_b | | \hat{\cal U}_b |}} \sum_L
\Xi_{\alpha}^L S^L_{ab} (\Xi_{\beta}^L)^* \ , \eqn\FSSform $$
where the sum is over the intersection of the two
untwisted stabilizers. The hats on $\hat{\cal S}$ and $\hat{\cal U}$
distinguish these groups from the (untwisted) stabilizers in the
${\cal H}$ extension, and we will use the symbol $\Xi$ to denote
$\hat{\cal U}$ characters. The ${\cal U}$ characters will be denoted
by $\Psi$. The labels $\alpha$ and $\beta$ distinguish the
resolved fixed points. For the resolved fixed points in the
${\cal H}$-extended theory we will use labels $i$ and $j$.

Let $N$ be the smallest integer larger than zero so that $J^N \in {\cal H}$.
Then $| {\cal M} | = N |{\cal H}| $. As was shown in the previous
section,
the field labelled $a^{\cal H}$ is a fixed point if and only if
there exists a current $h \in {\cal H}$ so that
$Jh a = a$, and $F(a,Jh,K) = 1$
for all $K\in {\cal U}_a$. If such an element $h$ exists we may
in fact choose a class representative $R_a(J)$ so that $R_a(J)a=a$ and
$F(a,R_a(J),K) = 1$ for all $K \in {\cal S}_a$. Clearly
$\hat{\cal S}_a$ consists of $N$ classes $(R_a(J))^n {\cal S}_a$,
$n=0,\ldots,N-1$, so that $|\hat{\cal S}_a|=N |{\cal S}_a|$.

Also the untwisted stabilizer ${\cal U}$
gets extended by $R_a(J)$ and all its powers, because
(omitting the arguments $a$ and $J$ for simplicity)
$$ F(R^nK,R^mK')=
F(R^n,R^m)
F({K},R^m)
F(R^n,K')
F(K,K')=F(K,K')\ , $$
for all $K$ and $K' \in {\cal S}_a$. Here we used the group property
in both arguments.
The result is equal to 1 if (and only if) in the second argument
$K'$ is restricted to ${\cal U}$. This shows that $R^m K'$ is
untwisted with respect to $\hat {\cal S}_a$ and hence is in $\hat {\cal U}_a$.

Hence in any case of interest to us $|\hat{\cal U}_a|=N |{\cal U}_a|$.

It should be noted that there are several other situations possible
that are not of immediate interest. For example, on some fields
the untwisted stabilizer may decrease rather than increase in size.
Then the current $J^{\cal H}$ acts non-trivially on the fixed point
resolution labels and re-combines several resolved fixed points into
a single field. The stabilizer ${\cal S}_a$ on the other hand can never
get smaller, but it can be enlarged by non-trivial divisors of $N$.
This happens if $J^{\cal H}$ does not fix $a^{\cal H}$, but some
power (smaller than $N$) does fix  $a^{\cal H}$. However, in
neither of these cases $a^{\cal H}$ is fixed by $J^{\cal H}$, and
hence we do not have to consider these possibilities.

Given \FSSform\ we can obtain the fixed point resolution
matrix $J^{\cal H}$ by ``pulling out" the $\Zbf_N$ characters
corresponding to the extension by $J^{\cal H}$. In general
characters of a discrete group do not simply ``factorize"
in terms of characters of a subgroup.
To do this correctly we need the formalism developed in Appendix C
to write the characters of $\hat{\cal U}_a$ as a product of
those of ${\cal U}_a$, those of $\Zbf_N$ and an additional phase factor.

\section{Determination of the fixed point resolution matrix}

Consider a
$\CFT^{\cal H}$-field $a^{\cal H}$.
Out of the coset belonging to $a^{\cal H}$ we choose one
representative $a$, which is a field in $\CFT$. This choice is
arbitrary. Now fix the representatives for all
$\CFT^{\cal H}$-currents that fix $a^{\cal H}$. This
set of currents forms
the full stabilizer ${\cal T}_{a^{\cal H}}^{\cal H}$.  The
representatives can be chosen
by decomposing  ${\cal T}_{a^{\cal H}}^{\cal H}$
into cyclic factors, and choosing representatives
for the
generator of each cyclic factor.
We can choose these
basis representatives $R_a(\hbox{basis})$
 so that they fix the field $a$ (rather than
fixing it up to an element of ${\cal H}$).
As was shown in Appendix B,
we can also choose
them so that the twist $F(a,R_a(\hbox{basis}),{\cal S}_a)=1$. Both properties
are preserved under multiplication of representatives.
It follows
then that the representatives $R_a(K)$ fix $a$ for all
$K^{\cal H} \in {\cal T}_{a^{\cal H}}^{\cal H}$, and that they are
untwisted w.r.t. the stabilizer ${\cal S}_a$.

Having fixed the representatives we
compute the phases $\phi$ that appear if we decompose
${\cal T}_a$ characters into ${\cal U}_a$ characters. Since
${\cal T}_a \supset \hat{\cal U}_a \supset {\cal U}_a$ this fixes
the decomposition of all possible $\hat{\cal U}_a$ groups that may
occur in various extensions of $\CFT^{\cal H}$.
When we extend by $J^{\cal H}$ we now have a canonical representative
$R_a(J)$ available for each field $a^{\cal H}$ (and representative $a$).
The phases $\phi$ depend on the choice of representatives $R_a$.

Using the results of appendix $C$ we can write
the characters of $\hat {\cal U}_a$ as follows
$$ \Xi^{R_a(J^n) h}_{m,i} =
\psi^n_m(\Zbf_N) \phi_i(n,a) \Psi^h_i({\cal U}_a) $$
where $h \in {\cal U}_a$ and $i$ labels the distinct
${\cal U}_a$ characters. To simplify the notation
we write $n$ instead of $J^n$ in the upper index of $\psi$
and the first argument of $\phi$.  In the following we simplify
the notation further by dropping the arguments $\Zbf_N$ and ${\cal U}_a$,
which are self-understood.

The formula for $S$ in the
${\cal M}$-extended theory is \FSSform, where $\alpha$ should
be interpreted as $(m,i)$.
We can split the sum into a sum over $J$-cosets and a sum over ${\cal H}$:
$$ S_{(a,\alpha)(b,\beta)}= { | {\cal M} | \over
\sqrt{ | \hat{\cal S}_a | | \hat{\cal U}_a |
| \hat{\cal S}_b | | \hat{\cal U}_b |}}\sum_{n=0}^{N-1}\ \left[ \ \
\sum_{h \in {\cal H}, J^n h  \in \hat{\cal U}_a \cap
\hat{\cal U}_b}
\Xi_{\alpha}^{J^n h} S^{J^n h}_{ab} (\Xi_{\beta}^{J^n h}
)^*\right] \ . $$

Clearly there can only be a contribution if it is possible to
choose $h \in {\cal H}$ in such a way that
$J^n h$ is in the intersection of $\hat{\cal U}_a$ and
$\hat{\cal U}_b$.  We will denote this particular
element $J^nh$ as $R_{ab}(J^n)$.
The complete set of such elements is
obtained by multiplying $R_{ab}(J^n)$ by any element of
${\cal U}_a \cap {\cal U}_b$. Hence we can now write the summation as
$$\eqalign{ S_{(a,\alpha)(b,\beta)}= { | {\cal M} | \over
\sqrt{ | \hat{\cal S}_a | | \hat{\cal U}_a |
| \hat{\cal S}_b | | \hat{\cal U}_b |}}  \sum_{n=0}^{N-1}\ \bigg[ \ \
\sum_{ h  \in {\cal U}_a \cap
{\cal U}_b}
&\Xi_{\alpha}^{R_{ab}(J^n) h}
 S^{R_{ab}(J^n) h}_{ab} \cr &(\Xi_{\beta}^{R_{ab}(J^n) h}
)^*\bigg]  \ .\cr} $$
A legitimate choice for $R_{ab}(J^n)$ is $R_{ab}(J^n)=[R_{ab}(J)]^n$.
However, the choice is irrelevant,
since legitimate choices of $R_{ab}(J^n)$
differ by elements of ${\cal U}_a \cap {\cal U}_b$, over which we sum.

The representatives $R_{ab}(J^n)$ are valid choices
for $\hat{\cal U}_a/{\cal U}_a$
as well as $\hat{\cal U}_b/{\cal U}_b$, but in general
not equal to the choices $R_a(J^n)$
and $R_b(J^n)$ made {\it a priori}. We cannot allow the choice
of representatives in $\hat{\cal U}_a/{\cal U}_a$ to depend on $b$.
We must
choose a fixed basis on each stabilizer, and refer everything
to that basis. In other words, we must replace $R_{ab}(J^n)$ by
either $R_a(J^n)$ or $R_b(J^n)$.
To do this we use
the group property
$$\Xi_{\alpha}^{R_{ab}(J^n)  h}
=\Xi_{\alpha}^{R_a(J^n) h}
    \Xi_{\alpha}^{(R_{ab}(J^n)/R_a(J^n))} \eqn\XiFac $$
and the same for the index $b$. Now it should be noted that
$R_{ab}(J^n)$ and $R_a(J^n)$ are both representatives of the
same coset class, and hence their ratio is an element of ${\cal H}$,
Furthermore, since both
representatives fix $a$
their ratio is in fact in ${\cal S}_a$. Finally, by construction
both $R_{ab}$ and $R_a$ are untwisted with respect to ${\cal S}_a$,
and therefore their ratio
$R_{ab}/R_a$ is an element
of ${\cal U}_a$.
Hence the second $\hat{\cal U}_a$
character $\Xi_{\alpha}$
in \XiFac\ reduces to a ${\cal U}_a$ character $\Psi_i$.

Now we wish to split off $\Zbf_N$ characters.
Writing $\alpha$ as $(m,i)$
$$\Xi_{\alpha}^{R_a(J^n) h}
    \Psi_{i}^{(R_{ab}(J^n)/R_a(J^n))}=
 \psi^n_m \phi_i(n,a) \Psi^h_i
\Psi_{i}^{(R_{ab}(J^n)/R_a(J^n))} $$
Substituting this into the general formula for $S$ we find
$$\eqalign{ S_{(a,m,i)(b,m',j)}=  { | {\cal M} | \over
\sqrt{ | \hat{\cal S}_a | | \hat{\cal U}_a |
| \hat{\cal S}_b | | \hat{\cal U}_b |}}   \sum_{n=0}^{N-1}\ \bigg[ \ \
&\sum_{ h  \in {\cal U}_a \cap
{\cal U}_b}S^{R_{ab}(J^n) h}_{ab} \cr
&\psi^n_m \phi_i(n,a)
\Psi_i^{(R_{ab}(J^n)/R_a(J^n))h}\cr
  &
[\psi^n_{m'} \phi_{j}(n,b)
\Psi_{j}^{(R_{ab}(J^n)/R_b(J^n))h}]^* \bigg]  \ .\cr} \eqn\finalform $$

{}From this expression we read off the fixed point matrix for
$J^n{\cal H}$, by removing the $\Zbf_N$ characters. The
formula \FSSform\ for the $\Zbf_N$ extension of $\CFT^{\cal H}$,
for matrix elements of untwisted fixed points of $J$ (\ie\ with
${\cal S}^{\cal H}_a = {\cal U}^{\cal H}_a = \Zbf_N$)
 reads
$$ S_{((a,i),m)((b,j),m')} = {N \over \sqrt{N^4}}
\sum_{n=0}^{N-1} \psi^n_m S^n_{(a,i)(b,j)} (\psi^n_{m'})^* $$
Comparing this with \finalform\ we read off the matrix elements
of $S^{J^{\cal H}}$ (note that the term
with $n=1$ in the sum
is the generator of the orbit,
which corresponds by construction to $J^{\cal H}$)
$$\eqalign{ S^{J^{\cal H}}_{(a,i)(b,j)}=
&{ | {\cal H} | \over \sqrt{
  | {\cal S}_a || {\cal U}_a |
| {\cal S}_b || {\cal U}_b | } }
\bigg[
\sum_{ h  \in {\cal U}_a \cap
{\cal U}_b}
\Psi^h_i S^{R_{ab}(J) h}_{ab}(\Psi_j^h)^*\bigg]  \cr
&  \phi_i(J,a)\Psi_i^{(R_{ab}(J)/R_a(J))}\ \
 [\phi_{j}(J,b)\Psi_{j}^{(R_{ab}(J)/R_b(J))} ]^*
\ .\cr} \eqn\MainForm $$
This comparison may seem to give us
the fixed point resolution matrices for
all powers of $J^{\cal H}$. This, however, is misleading because
if $N$ is not prime some powers of $J^{\cal H}$ may have
additional fixed points not seen here. Obviously
the fixed point resolution matrices for powers of $J^{\cal H}$
is given by the same formula, but with a larger index set
$(a,i)$.

Eqn. \MainForm\ is the main result of this paper. We will now examine
some consequences

\chapter{Properties of fixed point resolution matrices}

One can go systematically through conditions $\{1\}$ -- $\{6\}$
that fixed point resolution matrices must satisfy, but we will
omit most of the details here. Condition $\{1\}$ is satisfied
by construction, and checking conditions $\{2\}$ and $\{3\}$ is
mainly a tedious exercise, similar to the one presented in \FSSs, but
with a few extra twists.

Also checking $\{6\}$ is straightforward. One finds furthermore
that the matrices $S^J$ -- unlike those of unextended WZW models --
are in general not symmetric. An example is $A_{4,5}A_{4,5}$
extended with the current $(J,J)$, where
$J$ is a simple current of $A_{4,5}$.
 The new theory has a center
$\Zbf_5$ generated by (the orbit of) $(J,0)$. The fixed point
resolution matrices for $(J,0)$, $(J^2,0)$, $(J^3,0)$, $(J^4,0)$
are all different and all of them are asymmetric.

\section{Simple current twists}

Condition $\{4\}$ is of special interest, since it involves a new
quantity, the twist $F^{\cal H}$. The twist $F^{\cal H}(a^{\cal H},K^{\cal  
H},J^{\cal H})$
is unambiguously defined only if $K^{\cal H}$ fixes $a^{\cal H}$.
If we take for $K^{\cal H}$ a representative $R_a(K)$ that is untwisted
w.r.t. ${\cal U}_a$ and fixes $a$ we find immediately
$$ F^{\cal H}(a^{\cal H},K^{\cal H},J^{\cal H}) =
F(a,R_a(K),R_{ab}(J)h)=F(a,R_a(K),R_{ab}(J))$$
note that the result does not depend on $h$ (since $F(a,R_a(K),h)=1$) and
hence not on the choice of representative $R_{ab}(J)$.

However, it does appear that $F^{\cal H}$ now depends on both $a$ and $b$.
Although condition $\{5\}$ looks like merely a definition of $F$,
it does imply the non-trivial requirement that $F$ be independent
of $b$. In addition, for unextended WZW-models $F$ is also
independent of $a$, but this is not required.

Let us therefore carefully
examine the $a$ and $b$ dependence of $F^{\cal H}$.
Compare first $F(a,R_a(K),R_{ab}(J))$ and $F(a,R_a(K),R_{ab'}(J))$. Here we
made the induction hypothesis
 that $F$ has no {\it explicit} dependence on $b$.  The
only potential
$b$-dependence is then
via the choice of representatives $R_{ab}$. However, $R_{ab}$
and $R_{ab'}$ differ by an element of ${\cal U}_a$, and
since $F(a,R_a(K),{\cal U}_a)=1$ we have  
$F(a,R_a(K),R_{ab}(J))=F(a,R_a(K),R_{ab'}(J))$,
so that there is no $b$-dependence. Then the obvious choice is $b=a$,
so that (since $R_{aa}(J)=R_a(J)$)
$$ F^{\cal H}(a^{\cal H},K^{\cal H},J^{\cal H}) =
F(a,R_a(K),R_{aa}(J))=F(a,R_a(K),R_{a}(J))\ .$$

Now consider the $a$-dependence.
Even if $F$ does not depend on $a$ explicitly via its first argument,
there is no reason why an $a$-dependence could not arise
through the choice of representatives $R_a$.
Indeed, it is easy to find an example. Consider the
minimal $N=1$ super conformal field theory realized by the
coset CFT $SU(2)_4 \times SU(2)_2 / SU(2)_6$.
(the supersymmetry plays no r\^ole here, however).
The field
identification current is $(4,2,6)$, and there are three
non-trivial simple currents, whose orbits are
$J=(4,0,0)+(0,6,2)$, $K=(4,0,2)+(0,6,0)$ and $L=(0,0,2)+(4,6,0)$
(the latter
is in  fact the supercurrent). Consider the field representatives
$a=(2,1,1)$ and $a'=(4,3,1)$, which are both fixed by all three
currents. The current representatives that actually fix these
field representatives are respectively $R_a(J)=(4,0,0), R_a(K)=(4,0,2)$
and $R_a(L)=(0,0,2)$. For $a'$ we find
$R_{a'}(J)=(0,6,2)$, $R_{a'}(K)=(0,6,0)$ and $R_{a'}(L)=(0,0,2)$.
The twists all originate from the $SU(2)_k$ twists $F(a,J,J)=(-1)^{k/2}$
for $k$ even. It is then easy to see that all twists among
$J,K$ and $L$ are opposite for $a$ and $a'$.

\section{Conjugation twists}

To
determine the charge conjugation matrix
we need to compute $[S^{J^{\cal H}}]^2$,
from which the coefficients $\eta$ can be extracted.
After a lengthy computation we arrive at the following
complicated expression
$$ \eta^{J^{\cal H}}_{a^{\cal H}}  = \eta^{R_{ac}}_{K_a  
a}F(a,K_a,R_{ac})^*\phi_i(a)(\phi_{\pi^a(i)}(c)^*\
  \Psi^{R_{ac}/R_a}_i(\Psi^{R_{ac}/R_{c}}_{\pi^a(i)}  )^*\ .$$
The notation is as follows.
The conjugate of the orbit $a^{\cal H}$ is $c^{\cal H}$, where
$a$ and $c$ are arbitrarily chosen representatives. Since they are
arbitrarily chosen, they are in general not each other's conjugates,
but $c^*$ differs from $a$ by an element of ${\cal H}$. This
element is denoted as $K_a$, and its definition is
$K_a a = c^*$.
 Furthermore
$\pi^a$ is some -- in general $a$-dependent -- permutation,
defined by
$$F(a,K_a,h)^* \eta^{h}_{K_0a}
\Psi^{h}_i   = \psi^h_{\pi^a(i)} \eqn\pidef $$

This expression can be simplified considerably if we link the
choices of representatives on conjugate orbits.
If $c$ is conjugate to $a$, or lies on the same orbit as $a^*$,
a valid choice
is certainly
$$ R_{ac} =R_a=R_c \ . $$
Then automatically $\phi_i^a= \phi_i^c$.
Furthermore it is always possible to choose $\eta$ in such a
way that $\eta_{Ka}=\eta_a$.
Now $\eta^{\cal H}$ simplifies to
$$ \eta^{J^{\cal H}}_{a}  =  
\eta^{R_a(J)}_{a}F(a,K_a,R_a(J))^*\phi_i(a,J)\phi_{\pi^a(i)}(a,J)^*\
  $$
This expression may appear to depend on the choice of representatives
$R_a$ but actually this is not the case: all four factors depend
on $R_a$, but the overall dependence cancels.

The product formula for two $\eta$'s is
$$\eqalign{ \eta^{J^{\cal H}}_{a}\eta^{L^{\cal H}}_{a} &=
\eta^{R_a(J)}_{a}F(a,K_a,R_a(J))^*\phi_i(a,J)\phi_{\pi^a(i)}(a,J)^*\ \cr
&  \eta^{R_a(L)}_{a}F(a,K_a,R_a(L))^*\phi_i(a,L)\phi_{\pi^a(i)}(a,L)^*\  \cr
&=F(a,R_a(J),R_a(L))\eta^{R_a(J)R_a(L)}_a F(a,K_a,R_a(J)R_a(L)) \cr &
\phi_i(a,J)\phi_{\pi^a(i)}(a,J)^*\phi_i(a,L)\phi_{\pi^a(i)}(a,L)^* } $$
Here we used the formula
$$ \eta^J_a \eta^K_a = F(a,J,K) \eta^{JK}_a $$
as an induction hypothesis.

Now we may use $R_a(J)R_a(L)=R_a(JL) h(J,L)$. Since the set
of representatives was fixed
so that $R_a(M)$ fixes $a$ and is not twisted with respect
to ${\cal S}_a$ (where $M$ can be $K$,$L$ or $KL$), and was
determined up to elements in ${\cal U}_a$,
$h(J,L)$ must be an element of ${\cal U}_a$. Then
$$ \eta^{R_a(J)R_a(L)}_a = F(a,R_a(JL),h(J,L))\eta^{R_a(JL)}_a
\eta^{h(J,L)}_a = \eta^{R_a(JL)}_a
\eta^{h(J,L)}_a $$
and
$$ F(a,K_a,R_a(J)R_a(L))=F(a,K_a,R_a(JL))F(a,K_a,h(J,L)) $$
Furthermore
$$\phi_i(a,J)\phi_i(a,L)=\phi_i(a,JL)\Psi_i^{h(J,K)} $$
The $h(i,j)$-dependent phases all cancel because of \pidef,
and we are left with
$$ \eqalign{
 \eta^{J^{\cal H}}_{a}\eta^{L^{\cal H}}_{a}
&= F(a,R_a(J),R_a(L)) \eta^{(JL)^{\cal H}}_a\cr
&= F^{\cal H}(a^{\cal H},J^{\cal H},L^{\cal H}) \eta^{(JL)^{\cal H}}_a\cr}
 $$
This shows that the relation $G=F$ does indeed hold, if it was
valid in the unextended theory.

\chapter{Conclusions}

We have derived the generalization of the formula of \FSSs\ for
the modular transformation matrix of the zero-point holomorphic
blocks on the torus to modular transformation matrix of the one-point,
simple current blocks. In doing so,
we have also generalized \FSSs\ from unextended
WZW-models to extended ones. The main result, formula \MainForm,
gives the fixed point resolution matrices for all simple
current extended WZW-models and -- perhaps most importantly --
most coset CFT's.

The disclaimer made in \FSSs\ regarding these matrices applies here as
well: the conditions for $S^J$ are necessary, but not sufficient.
In particular fusion rule integrality does not
 -- or at least not manifestly --
follow from these conditions. Nevertheless, the fact that these matrices
provide natural and general solutions to the non-trivial conditions
of \FSSs\ may be regarded as strong evidence for their correctness.
Formula \MainForm\ has been checked numerically in many
examples, and gives integral fusion
rules in all cases that were examined.

In comparison to the unextended WZW-models
considered in \FSSs,
The matrices given \MainForm\ satisfy the conditions of \FSSs\ in
non-trivial ways, exploiting nearly all allowed features, with
the interesting exception of complex simple current twists.

\ack

I would like to thank the theory group at CERN,
where part of this work was done, for hospitality,
and J\"urgen Fuchs and Christoph Schweigert for their interest,
for many useful remarks and for carefully reading the manuscript.

\Appendix{A}

Here we show that if $F(a,K,J)$ is real, and
the fixed point resolution matrices of tensor products are
the tensor products of those of the factors, then $F=G$.

Suppose
two closed sets
of simple currents in CFT's ${\CFT}$ and ${\CFT'}$
and their twists on fields $a$ and $a'$ are isomorphic in the
following sense: their exists a one-to-one map $\pi(J)=J'$ from
the currents of ${\CFT}$ to the currents in ${\CFT'}$ which preserves
all fusion products, conformal weights modulo integers, and the
simple current twists,
\ie $F(a',\pi(J),\pi(K))=F(a,J,K)$. The currents are assumed
to have integral or half-integral spin, and to be local
with respect to each other. The group they form under fusion
will be called ${\cal H}$.

Denote fields in the tensor product $\CFT\otimes \CFT'$ as $(a,a')$.
In the tensor product we may   consider
the set of currents $(J, \pi(J))$. All currents of this form
are local with respect to each other, they have integral spin, they
form a closed set under fusion and on the field $(a,a')$ their
twists cancel. Hence ${\cal S}_{(a,a')} =  {\cal U}_{(a,a')} = {\cal H}$.
Then in order to satisfy $\{5b\}$ in the tensor product
theory we need
$$ G((a,a'),(K,K'),(L,L'))=G(a,K,L)G(a,K',L')=1$$
Hence knowing $G(a,K,L)$ for one theory, we have now derived
$G(a',K',L^{\prime})$ for any other, isomorphic set of currents. In particular,
if for any set of currents and twists we can find an example
(a CFT, a set of currents, and a field $a$)
with $F(a,K,L)=G(a,K,L)$ we have proved this relation for any CFT.

Since $G$ and $F$ depend on just two currents,
it is clearly sufficient to look only at pairs of currents,
\ie\ at $\Zbf_N \times \Zbf_M$. Let us call the generators of
the two factors $J$ and $K$.
Denote the spins modulo integers
by $s_J$ and $s_K$.
In general we have \FSSs\ $F(a,J,J)=(-1)^{2s_J}$. By the
properties of twists, all twists are then encoded in the quantity
$F(a,J,K)\equiv F$.
The following table lists the possibilities and
the realization of these possibilities in terms of a tensor product.
The factors in the tensor product are denoted
$A_{N-1}$ (for $A_{N-1}$ at level $N$) and ${\cal I}$
(for the Ising model). The first two lines describe theories
with a single cyclic factor, the remaining ones list all possibilities
for $\Zbf_N \times \Zbf_M$, assuming $F$ is real.
Note that $F(a,K,J)=-1$ requires $N$ (and $M$) to have even order,
since $F(a,K^N,J)=(-1)^N=1$.
For all these
theories \GisF\ does indeed hold.

\input tables
\vskip 1.truecm
\begintable
$s_J$ | $s_K$ | $N$ | $M$ | $F$ | factors
| $J$ | $K$ \cr
0 | -- | any | -- | -- | $A_{N-1}$ | $J$ | -- \nr
$\half$ | -- | even | -- | -- | $A_{N-1}\ {\cal I}$ | $(J,\Psi)$ | -- \nr
0 | 0 | any | any | 1 | $A_{N-1} A_{M-1}$ | $(J,1)$ | $(1,J)$
 \nr
0 | $\half$ | any | even | 1 | $A_{N-1} A_{M-1}  \ {\cal I}$
| $(J,1,1)$ | $(1,J,\Psi)$
 \nr
$\half$ | $\half$ | even | even | 1 | $A_{N-1} A_{M-1}  \ {\cal I}^2$|
$(J,1,\Psi,1)$ | $(1,J,1,\Psi)$
 \nr
0 | 0 | even | even | -1 | $A_{N-1} A_{M-1}\ {\cal I}^3 $ |
$(J,1,\Psi,\Psi,1)$ | $(1,J,\Psi,1,\Psi)$\nr
0 | $\half$ | even | even | -1 | $A_{N-1} A_{M-1}\ {\cal I}^2 $
|$(J,1,\Psi,\Psi)$ | $(1,J,\Psi,1)$\nr
$\half$ | $\half$ | even | even | -1 | $A_{N-1} A_{M-1}\ {\cal I} $ |
$(J,1,\Psi)$ | $(1,J,\Psi)$\endtable

Although most of this argument can be extended to complex $F$,
one cannot
rule out the possibility that $G(a,J,K)=1$ for
all $a$, $J$ and $K$ (which is obviously the only way out).
However, no examples of this kind are known.

\Appendix{B}

\leftline{\bf Untwisted representatives}

Here we will show that if the $\CFT^{\cal H}$
current $J^{\cal H}$ fixes a field $a^{\cal H}$, one can find an
orbit representative $R_a(J)$ in the ${\cal H}$-orbit of $J$ such that
$$ F(a,R_a(J),K)=1 \hbox{~for all~} K \in {\cal S}_a \ . $$

We   drop the explicit dependence
on $a$ henceforth. Consider some representative $X(J)$ that fixes $a$.
The choice of $X(J)$ is fixed up to
multiplication by elements of ${\cal S}$.
Consider ${\cal S}/{\cal U}$.
The twist $F(K,L)$ satisfies the group properties
$$ \eqalign{ F(K_1K_2,L)&=F(K_1,L)F(K_2,L)\cr
F(K,L_1L_2)&=F(K,L_1)F(K,L_2)\cr } $$
if $K,K_1$ and $K_2$ fix $a$ ($L,L_1$ and $L_2$ must also fix $a$, since
otherwise $F$ is not defined). Hence
$$ F(K{\cal U},L{\cal U}) = F(K,L) $$
\ie\ the twist is constant on coset classes.
Choose a set of generating representatives $L_i$ of
${\cal S}/{\cal U}$, so that
the full set of representatives is of the form
$L^{\vec m} \equiv \prod_i (L_i)^{m_i}$, $0 < m_i < N_i$, with
$N_i$ the smallest positive integer so that $L^{N_i} \in {\cal U}$.
Then it follows from the product rules of $F$ that the
relative twists of the generating representatives has the form
$$ F(L_i,L_j) = e^{2\pi i r_{ij} / N_{ij}}\ , $$
where $N_{ij}$ is the greatest common divisor of $N_i$ and $N_j$, and
$r_{ij}$ a set of integers.

The twist of the $X(J)$ with respect to the generators is
$$ F(X(J),L_i) = e^{2\pi i p_{i} / N_{iJ}}\ .$$
Our claim is now that another representative exists which has
trivial twist. Consider $X(J) L^{\vec k}$. Then
$$ F(X(J) L^{\vec k},L_i)=e^{2\pi i (p_i /N_{iJ} + k_j r_{ji}/N_{ij})} $$
Clearly we need to find a solution $\vec k$ to the equation
$$ \sum_j k_j {r_{ji}\over N_{ij}}= -{p_i \over N_{iJ}} \ \mod 1\eqn\KS $$
The definition of ${\cal S}/{\cal U}$ implies that
all currents $L^{\vec m}$ are twisted with respect at least on
other current. This means that $F(L^{\vec m},L_i)=1$ has no
solutions except the trivial one, $\vec m=0$. This implies that
$$ m_i {r_{ij} \over N_{ij}} = 0  \mod 1\ \ \Longleftrightarrow\ \ \vec m=0 \ .
 \eqn\vanish $$
Roughly speaking \KS\ has precisely one solution because \vanish\ implies
that the matrix $r_{ij}/N_{ij}$ is ``invertible". But because
we work with integers and modulo 1 this is not correct as it stands.

The correct solution to this problem can be found in \BeBU, where
exactly the same condition emerges in a quite different situation,
although also related to simple currents. In \BeBU\
the most general simple current invariant was constructed given
a center $\Zbf_{N_1}\ldots \Zbf_{N_n}$ with
generators $J_1\ldots J_n$. It was found that each
invariant was specified by a matrix
(called $X$ in the second paper in \BeBU) of the form $r_{ij}/N_{ij}$,
where again $r_{ij} \in \Zbf$ and $N_{ij}={\rm GCD}(N_i,N_j)$. In
that context, \vanish\ is the condition for having
a pure automorphism invariant without extension of the chiral algebra.
Equation \KS\ appears as the equation determining the way the
automorphism acts on a field of charges $p_i/N_{iJ}$, which are
allowed charges w.r.t. the currents $J_i$. The modular invariant
can be written down and checked explicitly, and since it is
of automorphism type there is one and only one solution $\vec k$
for any $\vec p$. This proves not only that the ``untwisted"
representative exists, but also that it is unique up
to multiplication by elements of ${\cal U}_a$.

\Appendix{C}

\leftline{\bf Subgroup Characters}

Consider a finite, discrete, abelian
group ${\cal G}$ with subgroup ${\cal H}$. The
coset ${\cal C}={\cal G}/{\cal H}$ is also an abelian group.
Both  ${\cal H}$ and ${\cal C}$ have a set of characters,
which we will denote by $\Psi$ and $\psi$ respectively.
We want to express the characters of ${\cal G}$ in
terms of $\Psi$ and $\psi$.

We denote the elements of ${\cal H}$ as $h$, and the character labels
as $i,j,\ldots$. The characters will be denoted as
$ \Psi^h_i$. Being characters, they satisfy the group properties
$$ \Psi_i^h \Psi_i^g=\Psi_i^{hg}\ , \ \ \ \Psi^1_i=1 \ . $$
Furthermore they satisfy an orthogonality relation
$$ \sum_h \Psi_i^h (\Psi_j^h)^* = |{\cal H}| \delta^{ij} $$
and a completeness relation
$$ \sum_i \Psi_i^h (\Psi_i^g)^* = |{\cal H}| \delta^{hg} \ . $$
The labels $i$ are not assumed to be ordered in any particular way.

The elements of ${\cal C}$ are denoted as $J$, and the
character labels are denoted $m,n,\ldots$. The characters are thus
denoted $ \psi^J_m$, and have properties analogous to those of $\Psi$.
In each  ${\cal G}/{\cal H}$ coset $J$ we choose a representative
$R(J) \in {\cal G}$.
The elements of ${\cal G}$ are
all of the form
$R(J)h$, for all $J \in {\cal C}$ and
all $h \in {\cal H}$.
In general
the set of representatives
only closes up to elements of ${\cal H}$:
$$ R(J)R(K)=R(JK) h(J,K)\ , \eqn\RepProd $$
with $h(J,K) \in {\cal H}$.

A complete set of character labels is
$(m,i)$ where $m$ and $i$ are from the same index sets as before. This
ensures that we get at least the right number of labels. An obvious
guess for the characters $\Xi$ of ${\cal G}$ could be to take simply
the product of those of ${\cal H}$ and ${\cal G}/{\cal H}$, but
it is easy to see that in general this yields incorrect group properties,
for any choice of representatives.
To remedy this we try the following ansatz
$$ \Xi^{R(J)h}_{(m,i)} = \psi^J_m \Psi^h_i \phi^J_i\ , $$
where $\phi$ is assumed to be a phase (it is actually a cocycle --
see the appendix of the third paper of \FuSc\ for a related discussion).

Let us first demonstrate orthogonality and completeness of these
characters.

\item\dash Orthogonality:\break
$$ \sum_{m,i}\Xi^{R(J)h}_{(m,i)}[\Xi^{R(J')h'}_{(m,i)}]^* =
\sum_{m,i}\psi^J_m \Psi^h_i \phi^J_i
[\psi^{J'}_m \Psi^{h'}_i \phi^{J'}_i ]^* $$
The sum over $m$ can be done and yields
(up to normalization)
$\delta_{JJ'}$ using orthogonality
in ${\cal C}$. Then the phases $\phi$ cancel, and we can perform the
sum over $i$.
\item\dash Completeness:\break
$$ \sum_{J,h}\Xi^{R(J)h}_{(m,i)}[\Xi^{R(J)h}_{(m',i')}]^* =
\sum_{J,h}\psi^J_m \Psi^h_i \phi^J_i
[\psi^{J}_{m'} \Psi^{h}_{i'} \phi^{J}_i ]^* $$
The idea is similar. Now we can sum
in ${\cal H}$ over $h$ to get $\delta_{ii'}$,
the phases $\phi$ cancel, and we then sum over $J$.

Note that it is essential in these arguments that $\phi$ does not depend on
on $h$ and $m$. Clearly orthogonality and completeness
do not fix $\phi$. Consider now the group property. On the
one hand
$$  \Xi^{R(J)h}_{(m,i)}\Xi^{R(K)g}_{(m,i)}=\Xi^{R(JK)h(J,K)hg}_{(m,i)}
=\psi^{JK}_m \Psi^{h(J,K)hg}_i \phi^{JK}_i$$
On the other hand the product is
$$ \psi^J_m \Psi^h_i \phi^J_i \psi^K_m \Psi^g_i \phi^K_i=
\psi^{JK}_m\Psi^{hg}_i\phi^J_i\phi^K_i $$
Comparing the two expressions we find
$$ \Psi^{h(J,K)}_i \phi^{JK}_i=\phi^J_i\phi^K_i \eqn\phidef $$
It is easy to show that a change in representatives in the cosets can be
compensated by a change in $\phi$.
This observation helps in determining the phases $\phi$, because
it allows us to make a convenient choice of representatives.
For the identity we choose representative $1$.
Choose a basis
 in ${\cal C}$ (by ``basis" we mean here a set of elements
$J_1, \ldots, J_n$ generating cyclic subgroups, so that any other element
can be written as
$ (J_1)^{m_1} \ldots (J_n)^{m_n} $).
Consider one such cyclic subgroup of order $N$,
 and denote the generator $J_{\ell}$
simply as $J$. If $R(J)$ is a representative of class $J$, then
$R(J)^m$ is a representative of class $J^m$. Within the cyclic subgroup
we have then
$$ R(J^m)R(J^n)=R(J^{m+n})$$
as long as $m+n < N$, the order of $J$. In general $R(J)^N = h_J \in {\cal H}$.
Then we find
$$ h(m,n)\equiv h(J^m,J^n)=1\ \ \hbox{if $m+n < N$};\ \
  h(J^m,J^n)=h_J\ \ \hbox{otherwise} $$
(note that $m+n < 2N$).
Now we can fix the representatives on all other elements as
$$ R((J_1)^{m_1} \ldots (J_n)^{m_n}) = R(J_1)^{m_1}\ldots R(J_n)^{m_n} $$
The computation of the function $h$ is then as follows
$$\eqalign{ R((J_1)^{m_1} \ldots (J_n)^{m_n}) &
R((J_1)^{k_1} \ldots (J_n)^{k_n})\cr&=R((J_1)^{m_1+k_1} \ldots (J_n)^{m_n+k_n})
  h_1(m_1,k_1) \ldots h_n(m_n,k_n)} $$

Now we determine the phases $\phi$. On one cyclic subgroup of order $N$,
and generated by $J$,
we have
$$ \eqalign{\phi^{J^2}_i &=\phi_i^J \phi_i^J \cr
 \phi^{J^3}_i &=\phi_i^{J^2} \phi_i^J =(\phi_i^{J})^3\cr
      &\ldots\cr
\phi_i^{J^{N-1}} &=\phi_i^{J^{N-2}}\phi_i^J=(\phi_i^{J})^{N-1} \cr
\phi_i^{J^{N}} &=\phi_i^{J^{N-1}}\phi_i^J [\Psi^{h_J}_i]^*=  
(\phi_i^{J})^{N}[\Psi^{h_J}_i]^* \cr}$$
On the other hand, since $J^N=1$ and the phases
$\phi$ depend only on the elements of ${\cal C}$, we must demand that
$$ \phi_i^{J^{N}} = 1 $$
This the implies
$$ \phi_i^J= [\Psi^{h_J}_i]^{1/N}\ , $$
i.e. it should be equal to one of the $N^{th}$ roots of
$\Psi^{h_J}_i$. It should not matter which root we choose, since
different choices amount only to permutations of labels $m$
on the ${\cal G}$ characters. Indeed, the different roots are
$$  [\Psi^{h_J}_i]^{1/N}\psi^J_m \ \  m = 0 \ldots, N-1\ . $$
With this definition the phases on each cyclic subgroup
satisfy
$$ \phi_i^{J^m}\phi_i^{J^n}=\phi_i^{J^{h+n}}\Psi^{h(m,n)}_i \ ,$$
with $h(m,n)$ as defined above. This procedure is repeated
for each cyclic subgroup generated by the basis elements $J_{\ell}$.
To keep the notation manageable we define the phases in
the ${\ell}^{\rm th}$ cyclic factor as
$$ \phi_i^m(\ell) \equiv \phi_i^{J^m_{\ell}}$$

Having in this manner determined the phases on the cyclic
subgroups of the individual generators, we have for the general case
$$ \phi_i^{m_1,\ldots,m_n} = \prod_{\ell=1}^n \phi^{m_{\ell}}_i({\ell}) $$
The final check is then the product rule
$$ \phi_i^{\vec m} \phi_i^{\vec n}
=\prod_{\ell=1}^n \phi^{m_{\ell}}_i({\ell}) \phi^{n_{\ell}}_i({\ell})
=\prod_{\ell=1}^n \phi_i^{m_{\ell}+n_{\ell}}
\Psi^{h_{\ell}(m_{\ell},n_{\ell})}_i
= \phi_i^{\vec m+\vec n} \Psi^{h(\vec m,\vec n)}_i$$
with the understanding that ${\vec m+\vec n}$ is taken modulo
$\vec N$.

Having determined $\phi$ for a convenient set of representatives $R$,
we can transform it to any other set in the following way.
Suppose we are given a any other set of representatives (with
the canonical choice on the identity)
$r(\vec m)$. Then $r(\vec m)=H(\vec m)R(\vec m)$,
with $H(\vec m) \in {\cal H}$, and
$R$ is the choice used above.
Then the phases $\phi$ with respect to the $r$-representatives are
$$ \phi_i^{\vec m}(r) = \phi_i^{\vec m}(R)
\Psi^{H(\vec m)}_i \ . $$
It is straightforward to check that these phases $\phi$ do indeed
satisfy \phidef, with $h(J,K)$ determined from $r(\vec m)$
using \RepProd.

Note that this character decomposition has the property that
it reduces to ${\cal H}$ characters for elements of ${\cal H}$, provided
we choose $R(1)=1$. Then $R(1)R(1)=R(1)$, hence $h(1,1)=1$ and
hence $\phi_i^1\phi_i^1=\phi^1_i$ so that $\phi^1_i=1$.

\par \penalty-4000\vskip\chapterskip
   \spacecheck\referenceminspace \immediate\closeout\referencewrite
   \referenceopenfalse
   \line{\fourteenrm\hfil REFERENCES\hfil}\vskip\headskip
   \endlinechar=-1
   \input referenc.texauxil
   \endlinechar=13
   
\end